\documentclass[11pt,a4paper]{article}
\usepackage[utf8]{inputenc}
\usepackage[T1]{fontenc}
\usepackage[margin=1in]{geometry}
\usepackage{amsmath,amssymb,amsthm}
\usepackage{mathrsfs}
\usepackage{authblk}
\usepackage[hidelinks]{hyperref}
\newtheorem{theorem}{Theorem}[section]
\newtheorem{lemma}[theorem]{Lemma}

\theoremstyle{definition}
\newtheorem{remark}[theorem]{Remark}

\newtheorem{definition}{Definition}[section]

\newcommand{\R}{\mathbb{R}}
\newcommand{\cL}{\mathcal{L}}
\newcommand{\cR}{\mathcal{R}}
\newcommand{\ov}[1]{\overline{#1}}

\newcommand{\W}{W_0^{s,p}(\Omega)}
\newcommand{\intO}{\int_\Omega}

\title{A rearrangement optimization problem for nonlocal $p$-Laplacian-type operators}

\author[1]{Chong Qiu}
\author[2]{Pengcheng Wu\thanks{Corresponding author.
  E-mail: \texttt{pcwu@suda.edu.cn} (Pengcheng Wu),
  \texttt{qchsuda@163.com} (Chong Qiu), \texttt{yuyingz@suda.edu.cn} (Yuying Zhou).}}
\author[2]{Yuying Zhou}
\affil[1]{School of Mathematics and Physics, Huai'an University, Huai'an 223003, China}
\affil[2]{Department of Mathematics, Soochow University, Suzhou 215006, China}
\date{}

\begin{document}

\maketitle

\begin{abstract}
We consider a rearrangement optimization problem related to the Dirichlet boundary problem driven by a general nonlocal operator $\cL_K$ of fractional
$p$-Laplacian type, with a nonnegative potential term, a nonlinear term and a source term
$f$ varying in the rearrangement class of a positive
function $f_0$. Under the condition that the potential term is dominated by the principal eigenvalue of the operator, we establish that the Dirichlet boundary problem has a unique solution by using the
Br\'ezis--Oswald and D\'iaz--Saa approaches, and we prove that the  associated energy functional admits a minimizer in the original rearrangement class so that we obtain a representation result of the
optimal source for the minimization problem.
\end{abstract}


\section{Introduction }\label{sec:intro}
Let $\Omega\subset\R^N$ ($N\ge 2$) be a bounded domain with $C^{1,1}$ boundary. Given a nonnegative function $f_0\in L^\infty(\Omega)$, we say that $f\in L^\infty(\Omega)$ lies in the rearrangement class of $f_0$, denoted $\cR$, if for all $t\ge 0$
\[
|\{x\in\Omega: f(x)>t\}| = |\{x\in\Omega: f_0(x)>t\}|,
\]
where $|\cdot|$ denotes the $N$-dimensional Lebesgue measure. Rearrangement classes provide a natural framework for optimization problems in the calculus of variations: one considers a functional $\Phi:\cR\to\R$ arising from a variational problem and studies the existence and characterization of extremal functions $\hat f\in\cR$ that maximize or minimize $\Phi$ over the class.

The systematic study of optimization problems on rearrangement classes was initiated by Burton in the seminal works \cite{B87,B89}, where the existence of maximizers for convex functionals was established through a weak$^*$ closure argument. A classical example is the following: for each $g\in\cR$ consider the Poisson problem
\[
\begin{cases}
-\Delta u = g(x) & \text{in } \Omega,\\
u = 0 & \text{on } \partial\Omega,
\end{cases}
\]
which admits a unique weak solution $u_g\in H^1_0(\Omega)$. The functional
\[
\Phi(g) = \int_\Omega g u_g\,dx
\]
was shown to admit a maximizer $\hat g\in\cR$, and moreover the optimal datum satisfies a pointwise relation $\hat g = \eta\circ u_{\hat g}$ for some nondecreasing function $\eta$, which has interesting consequences for the qualitative properties of the optimizer (see \cite{B87,B89,BurtonMcLeod91}).

Since these pioneering results, optimization problems in rearrangement classes have been investigated for a wide variety of elliptic equations. For instance, Cuccu, Porru, and Sakaguchi \cite{CuccuPorruSakaguchi11} studied the $p$-Laplacian equation with $p>1$ and a sublinear term:
\[
\begin{cases}
-\Delta_p u = g(x)u^{q-1} & \text{in } \Omega,\\
u = 0 & \text{on } \partial\Omega,
\end{cases}
\]
with $q\in[0,p)$, proving existence of both maximizers and minimizers of the functional $\Phi(g)=\int_\Omega g u_g^q\,dx$. In \cite{CuccuEmamizadehPorru09}, the minimization of the principal eigenvalue of the $p$-Laplacian with an indefinite weight was addressed. More results in this direction can be found in \cite{CuccuEmamizadehPorru06,EmamizadehPrajapat09,EmamizadehZivariRezapour07,Marras10,QiuHuangZhou15}.

In recent years, increasing attention has been devoted to nonlocal operators of fractional order (see \cite{DiNezzaPalatucciValdinoci12} for a comprehensive introduction). The fractional Laplacian $(-\Delta)^s$ ($s\in(0,1)$) and its nonlinear generalization, the fractional $p$-Laplacian
\begin{equation}\label{e01}
 (-\Delta)^s_p u(x) = C_{N,p,s}\lim_{\varepsilon\to 0^+}\int_{B_\varepsilon^c(x)}
\frac{|u(x)-u(y)|^{p-2}(u(x)-u(y))}{|x-y|^{N+ps}}\,dy,
\end{equation}
arise in diverse fields such as game theory, phase transitions, and anomalous diffusion (see \cite{BjorlandCaffarelliFigalli12,Caffarelli12,WuHuangZhou22}).

In the linear fractional setting, i.e., for $p=2$, rearrangement optimization problems have
been studied both for the principal eigenvalue and for the ground state energy. In
\cite{AneddaCuccuFrassu21}, the object of minimization is the first positive eigenvalue
$\lambda_1(\rho)$ of the weighted eigenvalue problem
\[
(-\Delta)^s u=\lambda\rho(x)u\ \text{ in }\Omega,\qquad u=0\ \text{ in }\Omega^c,
\]
where the weight $\rho\in L^\infty(\Omega)$ varies in the rearrangement class of a fixed
$\rho_0\in L^\infty(\Omega)$. As for the energy, Qiu, Huang and Zhou
\cite{QiuHuangZhou16} considered the Dirichlet problem
\[
-L^s_\theta u+h(x,u)=f(x)\ \text{ in }\Omega,\qquad u=0\ \text{ in }\Omega^c,
\]
where $L^s_\theta$ is the fractional Laplace-type operator
$L^s_\theta u(x)=\int_{\R^N}\frac{u(x+y)+u(x-y)-2u(x)}{|y|^{N+2s}}\theta(y)\,dy$
associated with a positive even weight $\theta$, and $h$ is nondecreasing in the second
variable with at most sublinear growth. Writing $H(x,t)=\int_0^th(x,\tau)\,d\tau$ and
denoting by $u_f$ the unique solution, they minimize the ground state energy
\[
\frac12\iint_{\R^N\times\R^N}\frac{(u_f(x)-u_f(y))^2}{|x-y|^{N+2s}}
\theta(x-y)\,dx\,dy+\int_\Omega H(x,u_f)\,dx-\int_\Omega fu_f\,dx
\]
as $f$ varies in the rearrangement class of a fixed $f_0\ge0$. Qiu and Zhou
\cite{QiuZhou2019} then treated the perturbed equation
\[
(-\Delta)^s u=f(x)+\lambda h(x)u\ \text{ in }\Omega,\qquad u=0\ \text{ in }\Omega^c,
\]
with $h>0$ a bounded weight, $f$ a given source and $0<\lambda<\lambda_1(h)$, where
$\lambda_1(h)$ is the first eigenvalue of $(-\Delta)^s u=\lambda h(x)u$; there the ground
state energy is the objective, minimized over the rearrangement classes of both $h$ and $f$.
Qiu \cite{Qiu2021} then treated the fractional $p$-Laplacian equation
\[
(-\Delta)^s_pu=f(x)+\lambda V(x)|u|^{p-2}u\ \text{ in }\Omega,\qquad u=0\ \text{ in }\Omega^c,
\]
with $p\ne2$, a fixed potential $V$ and a parameter $\lambda$, and showed that the associated
minimization problem over the rearrangement class of $f$ is solvable for a range of $\lambda$.

The nonlinear case was treated by Iannizzotto and Porru \cite{IannizzottoPorru2025}, who
considered nonlocal operators $\cL_K$ {defined in weak form by
\begin{equation}\label{eq:operator}
\langle\cL_K u,\phi\rangle
= \iint_{\R^N\times\R^N} |u(x)-u(y)|^{p-2}(u(x)-u(y))(\phi(x)-\phi(y))K(x,y)\,dx\,dy,
\end{equation}
for all $u,\phi\in W_0^{s,p}(\Omega)$, where $p>1$, $s\in(0,1)$, and $K:\R^N\times\R^N\to\R$ is a measurable kernel satisfying
\begin{itemize}
\item[(K1)] (symmetry) $K(x,y)=K(y,x)$ for a.e.\ $(x,y)\in\R^N\times\R^N$;
\item[(K2)] (two-sided bound) there exist constants $0<C_1\le C_2$ such that $C_1\le K(x,y)\,|x-y|^{N+ps}\le C_2$ for a.e.\ $(x,y)\in\R^N\times\R^N$.
\end{itemize}}
The first is the minimization of the principal eigenvalue
$\lambda_1(g)$ of the nonlinear eigenvalue problem
\[
\cL_K u=\lambda g(x)|u|^{p-2}u\ \text{ in }\Omega,\qquad u=0\ \text{ in }\Omega^c,
\]
as the weight $g$ varies in a rearrangement class. The second is the maximization of the
energy functional of the general Dirichlet problem
\[
\cL_K u+h(x,u)=g(x)\ \text{ in }\Omega,\qquad u=0\ \text{ in }\Omega^c,
\]
where $h$ is nondecreasing and grows at most like $|u|^{q-1}$ with $q<p$, and the
right-hand side $g$ varies in a rearrangement class.  {They proved the solvability of two
optimization problems  and obtained that  the optimizing datum
satisfies a pointwise relation $\hat g=\eta\circ u_{\hat g}$ with $\eta$ nondecreasing in \cite{IannizzottoPorru2025}.}  

 {Motivated by these works motioned above, we deals with a minimization of the energy functional for a 
equation of nonlocal operator $\cL_K$ of the type \eqref{eq:operator}  
 with a potential term and nonlinear term, as the source term varies in
a rearrangement class in this paper.

Consider  the following  Dirichlet problem involving  nonlocal operator $\cL_K$:}
\begin{equation}\label{eq:main}
\begin{cases}
\cL_K u - \lambda V(x)|u|^{p-2}u + h(x,u) = f(x) & \text{in } \Omega,\\[2mm]
u = 0 & \text{in } \Omega^c,
\end{cases}
\end{equation}
where $\lambda\ge 0$ is a parameter, $V\in L^\infty(\Omega)$ satisfies $V(x)\ge 0$ a.e.\ in
$\Omega$, $h:\Omega\times\R\to\R$ is a Carath\'eodory mapping, and $f$ belongs to the
rearrangement class $\cR$ of a fixed positive function $f_0$, i.e.\ $0<f_0\le M$ a.e.\ in
$\Omega$ for some $M>0$.
Denote  the energy functional as
\[
E(f,u) = \frac{1}{p}\Big( \iint_{\R^N\times\R^N} |u(x)-u(y)|^p K(x,y)\,dx\,dy - \lambda\int_\Omega V(x)|u|^p\,dx\Big)+\int_\Omega \Big[ H(x,u)-fu  \Big]\,dx ,
\]
where $H(x,t)=\int_0^t h(x,\tau)\,d\tau$. The potential term $-\lambda V(x)|u|^{p-2}u$
contributes the negative (destabilizing) summand $-\frac{\lambda}{p}\int_\Omega V|u|^p\,dx$
to the energy; under the spectral condition (V3) the $p$-homogeneous part of $E(f,\cdot)$ remains
positive definite. 
 For each $f\in\cR$, problem~\eqref{eq:main} admits a unique nonnegative weak solution $u_f\in W_0^{s,p}(\Omega)$ in the sense of Definition~\ref{def:weak} (Theorem~\ref{thm:exist}), and we set
\begin{equation*}
\Psi(f) = \inf_{u\in W_0^{s,p}(\Omega)} E(f,u) = E(f,u_f).
\end{equation*}
 In this paper we are interested in the following optimization problem:
\begin{equation}\label{eq:opt}
\min_{f \in \cR} \Psi(f) \tag{Opt}
\end{equation}

For the    nonlocal operator $\cL_K$, 
the assumptions (K1)--(K2) are exactly the kernel hypotheses of \cite{IannizzottoPorru2025}: no
translation invariance or radiality of $K$ is required, so that $\cL_K$ is allowed to be
anisotropic. The classical kernel $K_0(x,y)=|x-y|^{-N-ps}$ is the model example. 
When $C_1=C_2=C_{N,p,s}/2$, $\cL_K$ reduces exactly to the fractional $p$-Laplacian $(-\Delta)^s_p$ displayed in \eqref{e01}, which in turn coincides with the fractional Laplacian for $p=2$. The associated Gagliardo seminorm is
\[
[u]_K = \Big( \iint_{\R^N\times\R^N} |u(x)-u(y)|^p K(x,y)\,dx\,dy \Big)^{1/p}.
\]

We make the following hypotheses on $h$ and $V$ appearing in the Dirichlet problem \eqref{eq:main}.
\begin{itemize}
\item[(H1)] $h$ is a Carath\'eodory function, i.e., $x \mapsto h(x,t)$ is measurable for all $t \in \R$, and $t \mapsto h(x,t)$ is continuous for a.e. $x\in\Omega$.
\item[(H2)]
$h(x,\cdot)$ is nondecreasing on $\R$ for a.e.\ $x\in\Omega$, and $h(x,0)=0$ for a.e.\ $x\in\Omega$.
  \item[(H3)] There exist $C_0>0$ and $q\in(1,p_s^*)$ such that
  \[
  |h(x,t)|\le C_0(1+|t|^{q-1}) \quad \text{for a.e.\ } x\in\Omega,\ \forall t\in\R.
  \]
  \item[(H4)]  $t\mapsto h(x,t)/t^{p-1}$ is nondecreasing on $(0,\infty)$ for a.e.\ $x\in\Omega$.
\end{itemize}

\begin{itemize}
\item[(V1)] $\lambda\ge 0$.
\item[(V2)] $V\in L^\infty(\Omega)$ and $V(x)\ge 0$ for a.e.\ $x\in\Omega$.
\item[(V3)] $\lambda\|V\|_\infty < \lambda_1$, where $\lambda_1>0$ is the first eigenvalue of $\cL_K$ in $\Omega$ with Dirichlet boundary conditions:
\end{itemize}

\begin{equation}\label{eig}
\lambda_1 = \inf_{u\in W_0^{s,p}(\Omega)\setminus\{0\}}
\frac{[u]_K^p}{\int_\Omega |u|^p\,dx}.
\end{equation}

Two structural features distinguish \eqref{eq:main} from the equations studied in
\cite{IannizzottoPorru2025,QiuZhou2019}. First, 
the reaction term $h$ is a general nonlinearity subject to (H1)--(H4), whereas
\cite{QiuZhou2019} considers a linear perturbation $\lambda h(x)u$ whose weight varies in a
rearrangement class, and \cite{IannizzottoPorru2025} restricts to sublinear reactions. Second, the
potential term $-\lambda V(x)|u|^{p-2}u$, with $V$ fixed, is absent in
\cite{IannizzottoPorru2025}.
{When $p\neq 2$ the positivity of the $p$-homogeneous part of the energy does not imply that $E(f,\cdot)$ is convex. Hence the strict convexity argument of \cite{IannizzottoPorru2025} cannot be applied here.} We need to prove the uniqueness of the positive weak solution to (3) by new technology and method.
{A classical uniqueness criterion for semilinear elliptic equations is due to Br\'ezis and Oswald \cite{BrezisOswald}, who proved the uniqueness of positive solutions by using the monotonicity of the normalized reaction term.} D\'iaz and Sa\'a
\cite{DiazSaa} extended this approach to quasilinear elliptic
equations involving the $p$-Laplacian and established a comparison
inequality for positive solutions. In the nonlocal setting, Brasco and
Franzina \cite{BrascoFranzina2014} developed convexity principles and
Picone-type inequalities for local and nonlocal Dirichlet energies,
providing important tools for comparison arguments in nonlinear
nonlocal equations. Mugnai, Pinamonti, and Vecchi
\cite{MugnaiPinamontiVecchi} considered the standard fractional
$p$-Laplacian with nonlocal Robin boundary conditions and proved
uniqueness by truncating the solutions and adding $\varepsilon$ to the
denominators in order to construct regularized quotient test functions.
Kim \cite{Kim} used similar $\varepsilon$-regularized test functions
for a nonlocal Dirichlet problem with $p=2$ and applied the resulting
uniqueness theorem to a problem involving a discontinuous Kirchhoff
coefficient.
Unlike the above results, our setting allows both an arbitrary $p>1$
and a general symmetric measurable kernel $K(x,y)$; in particular, the
kernel is not required to be translation invariant or to have a radial
form. The equation also contains a negative $p$-homogeneous potential
term and a general nonlinearity $h$. We only assume that
$
t\longmapsto h(x,t)/t^{p-1}
$
is nondecreasing, rather than strictly increasing. The strictness needed
for uniqueness is supplied by the positive datum $f>0$. Our proof
uses one-sided $\varepsilon$-regularized test functions involving the
positive part, instead of the symmetric quotient test functions used in
the previous works. This construction requires neither a Hopf-type
boundary estimate nor a bound on the quotient of two solutions. Thus,
the Br\'ezis--Oswald uniqueness theory is extended here to
equations involving general measurable kernels, arbitrary $p>1$, a
negative potential term, and a reaction term satisfying a non-strict
monotonicity condition.

In order to obtain the results of the rearrangement optimization problem (Opt), we prove the absence of degenerate elements in  $ \ov{\cR}$ and the strong supporting hyperplane inequality, see Lemma \ref{lem:no-degenerate} and Lemma 
\ref{lem:hyperplane}. 
The resulting minimizers have the following structure (Theorem~\ref{thm:main}): the minimum of $\Psi$ in $\ov{\cR}$ 
is actually lied in the original class $\cR$ not in
$\ov{\cR}\setminus\cR$, and the optimal source $\hat f$ is related to the
corresponding solution by a pointwise relation $\hat f=\eta\circ u_{\hat f}$ with $\eta$
nondecreasing.

The paper is organized as follows. In Section~\ref{sec:prelim}, we introduce the necessary preliminaries on rearrangement classes and fractional Sobolev spaces. In Section~\ref{sec:exist},  by constructing  one-sided
$\varepsilon$-regularized comparison functions and exploiting  the
Br\'ezis--Oswald and D\'iaz--Saa approaches,  we deal with the existence and uniqueness of solutions to problem~\eqref{eq:main}.  In Section~\ref{sec:prelim-lemmas},  we first prove  the absence of degenerate elements in  $ \ov{\cR}$ and the strong supporting hyperplane inequality, then we study the solvability of the rearrangement optimization problem.

\section{Preliminaries}\label{sec:prelim}

In this section we collect some necessary preliminary results on rearrangement classes and on fractional Sobolev spaces, following the exposition of \cite{IannizzottoPorru2025}.

\subsection{Rearrangement classes}\label{subsec:rearrangement}

Let $\Omega\subset\R^N$ ($N\ge 2$) be a bounded domain, and let  $0< f_0\le M$ a.e.\ in $\Omega$ for some $M>0$. As defined in Section~\ref{sec:intro}, we say that $f\in L^\infty(\Omega)$ is a  rearrangement of $f_0$, denoted $f\sim f_0$, if for all $t\ge 0$
\[
|\{x\in\Omega: f(x)>t\}| = |\{x\in\Omega: f_0(x)>t\}|.
\]
The  rearrangement class of $f_0$ is
\[
\cR = \{f\in L^\infty(\Omega): f\sim f_0\}.
\]
{and its weak$^*$ closure in $L^\infty(\Omega)$ is denoted by} $\ov{\cR}$.
The following properties are proved in \cite[Theorem 6]{B87} and \cite[Lemma 2.2]{B89}:
\begin{itemize}
\item $\ov{\cR}$ is sequentially weak$^*$ compact;
\item $\ov{\cR}$ is convex;
\item every $f\in\ov{\cR}$ satisfies $0\le f\le M$ a.e.\ in $\Omega$.
\end{itemize}

Therefore, given a sequentially weak$^*$ continuous functional $\Gamma:\ov{\cR}\to\R$, there exist $\check f,\hat f\in\ov{\cR}$ such that for all $f\in\ov{\cR}$
\[
\Gamma(\check f) \le \Gamma(f) \le \Gamma(\hat f).
\]
In general, the extrema are not attained at points of $\cR$; a significant part of our analysis is devoted to proving that, for the specific functional $\Psi$, the minimizer actually lies in $\cR$.

About the optimization of linear functionals over $\ov{\cR}$, we need the following result, which also provides a representation formula for the optimizer:

\begin{lemma}\label{lem:linear-opt}
Let $u\in L^1(\Omega)$. Then:
\begin{itemize}
\item[\rm{(i)}] there exists $\hat f\in\cR$ such that $\displaystyle\int_\Omega \hat f u\,dx \ge \int_\Omega f u\,dx$ for all $f\in\ov{\cR}$;
\item[\rm{(ii)}] if $\hat f$ is unique, then there exists a nondecreasing function $\eta:\R\to\R$ such that $\hat f = \eta\circ u$ a.e.\ in $\Omega$.
\end{itemize}
\end{lemma}

This lemma follows from Burton's results \cite{B87,B89} and will be used in the proof of the main theorem.

\subsection{Fractional Sobolev spaces and the operator $\cL_K$}\label{subsec:sobolev}

For $p>1$, $s\in(0,1)$, and any open set $\Omega\subseteq\R^N$, we define the Gagliardo seminorm
\[
[u]_{s,p,\Omega} = \Big( \iint_{\Omega\times\Omega} \frac{|u(x)-u(y)|^p}{|x-y|^{N+ps}}\,dx\,dy \Big)^{1/p}.
\]
The fractional Sobolev space is defined by
\[
W^{s,p}(\Omega) = \{ u\in L^p(\Omega) : [u]_{s,p,\Omega} < \infty \}.
\]

When $\Omega$ is bounded with $C^{1,1}$-smooth boundary, we incorporate the Dirichlet condition by setting
\[
W_0^{s,p}(\Omega) = \{ u\in W^{s,p}(\R^N) : u=0 \text{ in } \Omega^c \},
\]
endowed with the norm $\|u\|_{W_0^{s,p}(\Omega)} = [u]_{s,p,\R^N}$. This is a uniformly convex, separable Banach space with dual $W^{-s,p'}(\Omega)$, and $C_c^\infty(\Omega)$ is dense in $W_0^{s,p}(\Omega)$. For the embeddings of $W_0^{s,p}(\Omega)$ into Lebesgue spaces, set
\[
p_s^* = \begin{cases}
\frac{Np}{N-ps} & \text{if } ps < N,\\
\infty & \text{if } ps \ge N,
\end{cases}
\]
and note that the value ``$\infty$'' has to be read with care: if $ps<N$ the embedding
$W_0^{s,p}(\Omega)\hookrightarrow L^q(\Omega)$ is continuous for all $q\in[1,p_s^*]$ and
compact for all $q\in[1,p_s^*)$; if $ps=N$ it is continuous and compact for every finite
$q\ge1$, but $W_0^{s,p}(\Omega)$ is not embedded in $L^\infty(\Omega)$; if $ps>N$ it
is continuous and compact for every $q\in[1,\infty]$, the endpoint $q=\infty$ following
from the embedding of $W_0^{s,p}(\Omega)$ into a H\"older space.
For a detailed account on fractional Sobolev spaces, we refer to \cite{DiNezzaPalatucciValdinoci12}.

Now let $K:\R^N\times\R^N\to\R$ be a measurable kernel satisfying (K1)--(K2). We introduce an equivalent norm on $W_0^{s,p}(\Omega)$ by setting
\[
[u]_K = \Big( \iint_{\R^N\times\R^N} |u(x)-u(y)|^p K(x,y)\,dx\,dy \Big)^{1/p}.
\]
Condition (K2) implies that $[\,\cdot\,]_K$ is equivalent to the standard Gagliardo seminorm $[\,\cdot\,]_{s,p,\R^N}$.

We now make precise the notion of solution used throughout the paper.

\begin{definition}[Weak solution]\label{def:weak}
Let (K1)--(K2) and (V1)--(V2) hold and let $f\in L^\infty(\Omega)$. A function
$u\in\W$ is called a weak solution of problem~\eqref{eq:main} if
\begin{equation}\label{eq:weak}
\langle\cL_K u,\varphi\rangle
-\lambda\int_\Omega V(x)|u|^{p-2}u\,\varphi\,dx
+\int_\Omega h(x,u)\,\varphi\,dx
=\int_\Omega f\,\varphi\,dx
\end{equation}
for every test function $\varphi\in\W$, where $\langle\cL_Ku,\varphi\rangle$ is given
by~\eqref{eq:operator}.
\end{definition}

\begin{remark}\label{rem:weak-well-posed}
Under the additional assumptions (H1)--(H3) all the terms occurring
in~\eqref{eq:weak} are finite for $u,\varphi\in\W$: the first one by H\"older's
inequality, which gives
$|\langle\cL_Ku,\varphi\rangle|\le[u]_K^{p-1}[\varphi]_K$; the second and the fourth
since $V,f\in L^\infty(\Omega)$, $|\Omega|<\infty$ and $u,\varphi\in L^p(\Omega)
\subset L^1(\Omega)$; and the third by the growth condition (H3) together with the
continuous embedding $\W\hookrightarrow L^q(\Omega)$, which yields
$|h(\cdot,u)|\in L^{q/(q-1)}(\Omega)$ and hence $h(\cdot,u)\varphi\in L^1(\Omega)$.
Moreover, for $u\in\W$ fixed,~\eqref{eq:weak} is precisely the Euler--Lagrange
equation of the energy functional $E(f,\cdot)$ displayed in~\eqref{eq:energy} below, that is,
$\partial_u E(f,u)=0$ in $W^{-s,p'}(\Omega)$; this equivalence is used in the proof of
Theorem~\ref{thm:exist}. Finally, since $\W$ is a vector space
and~\eqref{eq:weak} is required for every $\varphi\in\W$, testing against
$\varphi\in\W\cap L^\infty(\Omega)$ only would give the same notion.
\end{remark}

\section{Existence and uniqueness for the Dirichlet problem \eqref{eq:main}}\label{sec:exist}

In this section we prove the following existence and uniqueness result for the Dirichlet problem \eqref{eq:main}.

\begin{theorem}\label{thm:exist}
Let \textup{(K1)--(K2)}, \textup{(H1)--(H4)} and \textup{(V1)--(V3)} hold. Then for every $f$ with $0<f\le M$ a.e.\ in $\Omega$ for some $M>0$, problem~\eqref{eq:main} admits a unique nonnegative weak solution $u_f\in W_0^{s,p}(\Omega)$ and $E(f,u_f)=\inf_{u\in W_0^{s,p}(\Omega)} E(f,u)$.
\end{theorem}

Define the $p$-homogeneous quantity
\begin{equation*}
Q(u) = [u]_K^p - \lambda\int_\Omega V(x)|u|^p\,dx,
\end{equation*}
so that the energy functional takes the form $E(f,u)=\frac1pQ(u)+\int_\Omega[H(x,u)-fu]\,dx$.

In order to prove Theorem \ref{thm:exist}, we prove some lemmas first.

\begin{lemma}\label{lem:phi-bound}
Under assumptions \textup{(K1)--(K2)} and \textup{(V1)--(V3)},
\[
\delta [u]_K\le Q(u)^{1/p}\le [u]_K, \quad\forall u\in W_0^{s,p}(\Omega),
\]
where $\delta = (1 - \frac{\lambda\|V\|_\infty}{\lambda_1})^{1/p} > 0$.
\end{lemma}

\begin{proof}
For the upper bound, since $\lambda\ge 0$ and $V(x)\ge 0$ a.e.\ in $\Omega$, we have
\[
Q(u) = [u]_K^p - \lambda\int_\Omega V(x)|u|^p\,dx \le [u]_K^p,
\]
hence $Q(u)^{1/p}\le[u]_K$ for all $u\in W_0^{s,p}(\Omega)$.

For the lower bound, we get from (V3) that
\begin{align*}
Q(u) &\ge [u]_K^p - \lambda\|V\|_\infty\int_\Omega |u|^p\,dx\\
&\ge [u]_K^p - \frac{\lambda\|V\|_\infty}{\lambda_1}[u]_K^p
= \left(1 - \frac{\lambda\|V\|_\infty}{\lambda_1}\right)[u]_K^p = \delta^p[u]_K^p.
\end{align*}
This completes the proof.
\end{proof}

In what follows, $\|\cdot\|_p$ always denotes the $L^p(\Omega)$ norm.

\begin{lemma}\label{lem:ge0}
Under assumptions \textup{(H1)--(H2)}, we have $H(x,t)= \int_0^t h(x,\tau)\,d\tau\ge 0$.
\end{lemma}
\begin{proof}
  By (H1)--(H2), if $t>0$, then $h(x,t)\ge h(x,0)=0$ a.e.\ $x\in \Omega$, hence $H(x,t)=\int_0^t h(x,\tau)\,d\tau\ge 0$. If $t<0$, then $h(x,t)\le0$ a.e.\ $x\in \Omega$, which again yields $H(x,t)\ge 0$.
\end{proof}
\begin{lemma}\label{functional}
Under conditions \textup{(K1)--(K2)}, there exist $\alpha>1$ and $S>0$ such that for each $u\in\W$
\begin{equation}\label{sobolev}
 \|u\|_{p\alpha}^p\le S[u]_K^p.
\end{equation}
Moreover, the infimum in \eqref{eig} is positive and attained, so that
\begin{equation}\label{poincare}
 \lambda_1\|u\|_p^p\le[u]_K^p.
\end{equation}
\end{lemma}
The conclusions are direct consequences of the embedding
inequalities and Lemma \ref{lem:phi-bound}. We state them here for convenience in the
subsequent proofs.

\begin{lemma}\label{composition}
Let $u,v\in\W\cap L^\infty(\Omega)$ be nonnegative and let $L\ge\max\{\|u\|_\infty,\|v\|_\infty\}$. If $F:[0,L]^2\to\R$ is Lipschitz and $F(0,0)=0$, then $F(u,v)\in\W$.
\end{lemma}
This is the classical Lipschitz composition property of fractional Sobolev spaces for $0<s<1$ (see \cite{RS} for the discussion of composition operators); combined with $F(0,0)=0$ it yields the zero exterior condition. The two-variable version needed here also follows directly from
$|F(a,b)-F(c,d)|\le C(|a-c|+|b-d|)$ and the definition of the Gagliardo seminorm. Similarly, a globally Lipschitz function of one variable with $F(0)=0$ preserves $\W$.

\begin{lemma}\label{bounded}
Let \textup{(K1)--(K2)}, \textup{(H1)--(H2)} and \textup{(V1)--(V3)} hold, and let $f\in L^\infty(\Omega)$
satisfy $0\le f\le M$ a.e.\ in $\Omega$ for some $M>0$. Then every weak solution
$u\in\W$ of problem~\eqref{eq:main}, in the sense of Definition~\textup{\ref{def:weak}}, is
nonnegative and belongs to $L^\infty(\Omega)$.
\end{lemma}
\begin{proof}
Let $u\in\W$ be a weak solution
in the sense of~\eqref{eq:weak} and put $u^-=\min\{u,0\}$. Take $\phi=u^-\in\W$ in~\eqref{eq:weak}. Write $J_p(t)=|t|^{p-2}t$. For arbitrary real numbers $a,b$, a case distinction according to their signs gives
\[
 J_p(a-b)(a^--b^-)\ge|a^--b^-|^p.
\]
By (H2), $h(x,u)u^-\ge0$, while, since $f\ge0$ a.e.\ and $u^-\le0$, the right-hand side is $\intO fu^-\le0$. Hence
\[
 [u^-]_K^p-\lambda\intO V|u^-|^p\le0.
\]
Combining Lemma \ref{lem:phi-bound} yields $u^-=0$.

We now prove boundedness in detail. Set $a_V=\lambda\|V\|_\infty$. By (H2), $u\ge0$ implies $h(x,u)\ge0$. Thus for every nonnegative $\phi\in\W$ the weak equation gives
\begin{equation}\label{subsolution}
 \begin{aligned}
 \langle\mathcal L_Ku,\phi\rangle
 &=\intO\bigl(\lambda V(x)u^{p-1}+f(x)-h(x,u)\bigr)\phi\,dx\\
 &\le\intO(a_Vu^{p-1}+M)\phi\,dx.
 \end{aligned}
\end{equation}
Here the last step uses $\phi\ge0$, $h(x,u)\ge0$, $V\le\|V\|_\infty$ a.e.\ in
$\Omega$, and $0\le f\le M$ a.e.\ in $\Omega$.

\medskip\noindent
\textit{Step 1: Fixing the Sobolev exponent.}
By Lemma \ref{functional}, fix $\alpha>1$ and $S>0$ such that
\[
 \|w\|_{p\alpha}^p\le S[w]_K^p,\qquad w\in\W.
\]
When $sp<N$ one may take $1<\alpha<p_s^*/p$; when $sp\ge N$ one may take any fixed finite $\alpha>1$. 

\medskip\noindent
\textit{Step 2: Constructing admissible truncated test functions.}
Fix $r\ge p$ and set $\beta=r-p\ge0$. Given $\ell>0$, let
\begin{equation}\label{truncfunctions}
 \begin{gathered}
 u_\ell=\min\{u,\ell\},\\
 F_\ell(t)=t\min\{t,\ell\}^{\beta},\qquad
 G_\ell(t)=t\min\{t,\ell\}^{\beta/p},\qquad t\ge0.
 \end{gathered}
\end{equation}
When $\beta=0$ we adopt the convention $F_\ell(t)=G_\ell(t)=t$. These two functions are chosen because, once the truncation is removed,
\[
 F_\ell(u)\longrightarrow u^{r-p+1},\qquad
 G_\ell(u)\longrightarrow u^{r/p}.
\]
At this stage one cannot test the equation directly with $u^{r-p+1}$, since it has not yet been shown to belong to the energy space.

For fixed $\ell,\beta$, the functions $F_\ell,G_\ell$ are Lipschitz on $[0,\infty)$ and vanish at the origin. Defining them to be zero on the negative half-line yields globally Lipschitz functions. By the one-variable composition property following Lemma \ref{composition},
\[
 F_\ell(u),G_\ell(u)\in\W.
\]
Hence the nonnegative function $F_\ell(u)=uu_\ell^\beta$ is an admissible test function. The Lipschitz constants here may depend on $\ell,\beta$, since admissibility is only required for each fixed truncation.

\medskip\noindent
\textit{Step 3: Relating the differences of the two truncated functions.}
A direct computation gives, for almost every $t>0$,
\[
 \begin{array}{c|cc}
 &0<t<\ell&t>\ell\\ \hline
 F_\ell'(t)&(\beta+1)t^\beta&\ell^\beta\\
 G_\ell'(t)&(1+\beta/p)t^{\beta/p}&\ell^{\beta/p}.
 \end{array}
\]
Therefore
\begin{equation}\label{derivativebound}
 (G_\ell'(t))^p\le(1+\beta/p)^pF_\ell'(t).
\end{equation}
Take arbitrary $A>B\ge0$. By H\"older's inequality and \eqref{derivativebound},
\begin{align*}
 |G_\ell(A)-G_\ell(B)|^p
 &=\left|\int_B^A G_\ell'(t)\,dt\right|^p\\
 &\le(A-B)^{p-1}\int_B^A(G_\ell'(t))^p\,dt\\
 &\le(1+\beta/p)^p(A-B)^{p-1}
       \bigl(F_\ell(A)-F_\ell(B)\bigr).
\end{align*}
The same holds after interchanging $A$ and $B$; for $A=B$ both sides vanish. Hence for all $A,B\ge0$,
\begin{equation}\label{differencebound}
 |G_\ell(A)-G_\ell(B)|^p
 \le(1+\beta/p)^pJ_p(A-B)\bigl(F_\ell(A)-F_\ell(B)\bigr).
\end{equation}
Taking $A=u(x)$, $B=u(y)$, multiplying by $K(x,y)\ge0$ and integrating, we obtain
\begin{equation}\label{energybound}
 [G_\ell(u)]_K^p
 \le(1+\beta/p)^p\langle\mathcal L_Ku,F_\ell(u)\rangle.
\end{equation}
This step links the energy of $G_\ell(u)$ with the admissible test function $F_\ell(u)$ appearing in the equation.

\medskip\noindent
\textit{Step 4: Estimating the right-hand side by means of the equation.}
Taking $\phi=F_\ell(u)$ in \eqref{subsolution} and combining with \eqref{energybound} gives
\[
 [G_\ell(u)]_K^p
 \le(1+\beta/p)^p
 \intO(a_Vu^{p-1}+M)uu_\ell^\beta\,dx.
\]
Assume now that $u\in L^r(\Omega)$ is known. Since $u_\ell\le u$ and $r=p+\beta$,
\[
 u^pu_\ell^\beta\le u^r,\qquad
 uu_\ell^\beta\le u^{1+\beta}\le1+u^r.
\]
The last inequality can be checked separately on $0\le u\le1$ and on $u>1$: on the former it uses $u^{1+\beta}\le1$, on the latter $1+\beta\le r$. Therefore
\[
 \intO(a_Vu^{p-1}+M)uu_\ell^\beta\,dx
 \le(a_V+M)\|u\|_r^r+M|\Omega|.
\]
Using $1+\beta/p=r/p$ and the Sobolev inequality, we obtain
\begin{equation}\label{truncatedestimate}
 \|G_\ell(u)\|_{p\alpha}^p
 \le S(r/p)^p\bigl((a_V+M)\|u\|_r^r+M|\Omega|\bigr).
\end{equation}
In particular, the right-hand side is independent of the truncation level $\ell$.

\medskip\noindent
\textit{Step 5: Removing the truncation and increasing the integrability exponent.}
Letting $\ell\to\infty$ we have $G_\ell(u)\uparrow u^{r/p}$. Applying the monotone convergence theorem \cite[Chapter 1]{Rudin} to $G_\ell(u)^{p\alpha}$, estimate \eqref{truncatedestimate} gives
\[
 \left(\intO u^{\alpha r}\,dx\right)^{1/\alpha}
 \le S(r/p)^p\bigl((a_V+M)\|u\|_r^r+M|\Omega|\bigr).
\]
The argument allows the left-hand side to be infinite a priori; since the right-hand side is finite, the left-hand side is finite as well, whence $u\in L^{\alpha r}(\Omega)$. Moreover, the left-hand side is exactly $\|u\|_{\alpha r}^r$, so there exists $C>0$, independent of $r$ and $\ell$, such that
\begin{equation}\label{iteration}
 \|u\|_{\alpha r}^r
 \le Cr^p\bigl(\|u\|_r^r+|\Omega|\bigr).
\end{equation}
For instance one may take $C=Sp^{-p}(a_V+M)$, since $M\le a_V+M$. This yields simultaneously the increase of integrability and the corresponding norm control.

\medskip\noindent
\textit{Step 6: Iteration, and a uniform bound for the norms.}
Let
\[
 r_j=p\alpha^j,\qquad A_j=\max\{1,\|u\|_{r_j}\},
 \qquad j=0,1,2,\ldots.
\]
Since $u\in\W\subset L^p(\Omega)$, starting from $r_0=p$ and applying \eqref{iteration} successively gives $u\in L^{r_j}(\Omega)$ for all $j$. However, these integrability statements alone do not imply boundedness; the norms must be controlled.

Since $A_j\ge1$,
\[
 \|u\|_{r_j}^{r_j}+|\Omega|
 \le(1+|\Omega|)A_j^{r_j}.
\]
Substituting into \eqref{iteration} and taking $r_j$-th roots gives
\[
 A_{j+1}\le
 \max\{1,[C(1+|\Omega|)r_j^p]^{1/r_j}\}A_j.
\]
We now unroll the recursion step by step. Set
\[
 D_0=\max\{1,C(1+|\Omega|)\}.
\]
Since $D_0\ge1$ and $r_j>1$ we have $(D_0r_j^p)^{1/r_j}\ge1$; and since $C(1+|\Omega|)\le D_0$,
\[
 \max\{1,[C(1+|\Omega|)r_j^p]^{1/r_j}\}
 \le(D_0r_j^p)^{1/r_j}.
\]
Hence
\begin{equation}\label{recursionexplicit}
 A_{j+1}\le(D_0r_j^p)^{1/r_j}A_j.
\end{equation}
This implies that, for every $j\ge1$,
\begin{equation}\label{productestimate}
 A_j\le A_0\prod_{i=0}^{j-1}(D_0r_i^p)^{1/r_i}=A_0\exp\left(
 \sum_{i=0}^{j-1}\frac{\log D_0+p\log r_i}{r_i}
 \right).
\end{equation}
Moreover, since $r_i=p\alpha^i$ with $\alpha>1$,
$$
L_0:= \sum_{i=0}^{j-1}\frac{\log D_0+p\log r_i}{r_i}<\infty.
$$
Substituting back into
\eqref{productestimate}, we obtain
\[
 A_j\le A_0e^{L_0}\quad\text{for all }j\ge0.
\]
Setting $C_*=e^{L_0}$ and $C_\infty=C_*\max\{1,\|u\|_p\}$, we arrive at
\begin{equation}\label{uniformnorm}
 \|u\|_{r_j}\le A_j\le C_\infty\quad\text{for all }j.
\end{equation}
The constant $C_*$ depends only on the fixed Sobolev constant, on $p,\alpha,a_V,M$ and $|\Omega|$, and not on $j$ or on the truncation level. In other words, although each iteration step multiplies by a factor, the accumulated product of these factors still has a uniform finite bound, so that the $L^{r_j}$ norms do not grow without bound as $j$ increases.

\medskip\noindent
\textit{Step 7: Essential boundedness from the uniform norm control.}
Take any $T>C_\infty$ and set $\Omega_T=\{x\in\Omega:u(x)>T\}$. By \eqref{uniformnorm},
\[
 T^{r_j}|\Omega_T|
 \le\int_{\Omega_T}u^{r_j}\,dx
 \le\|u\|_{r_j}^{r_j}\le C_\infty^{r_j}.
\]
Thus
\[
 |\Omega_T|\le(C_\infty/T)^{r_j}\longrightarrow0\qquad(j\to\infty),
\]
since $C_\infty/T<1$ and $r_j\to\infty$. Hence $|\Omega_T|=0$. Taking $T=C_\infty+1/n$, $n=1,2,\ldots$, and taking the countable union of these null sets gives $u\le C_\infty$ almost everywhere. Therefore
\[
  u\in L^\infty(\Omega),\qquad
 \|u\|_\infty\le C_*\max\{1,\|u\|_p\}.
\]
\end{proof}

\begin{lemma}[Discrete Picone inequality]\label{picone}
{Writing $J_p(t)=|t|^{p-2}t$, we have, for all $a,b>0$ and $c,d\ge0$,}
\begin{equation}\label{piconeeq}
 J_p(a-b)\left(\frac{c^p}{a^{p-1}}-\frac{d^p}{b^{p-1}}\right)\le|c-d|^p.
\end{equation}
\end{lemma}
This is quoted directly from \cite[Proposition 4.2]{BrascoFranzina2014}, with the exponent there taken as $q=p$.

\begin{lemma}[D\'iaz--Saa type comparison]\label{ds}
Let $u,v\in\W\cap L^\infty(\Omega)$ be nonnegative. For $\varepsilon>0$ write, on all of $\R^N$,
\[
 u_\varepsilon=u+\varepsilon,\quad v_\varepsilon=v+\varepsilon,\quad
 z=(u_\varepsilon^p-v_\varepsilon^p)^+,\qquad
 \phi_u=\frac{z}{u_\varepsilon^{p-1}},\quad\phi_v=\frac{z}{v_\varepsilon^{p-1}}.
\]
Then $\phi_u,\phi_v\in\W$, and
\begin{equation}\label{dsineq}
 \langle\mathcal L_Ku,\phi_u\rangle-
 \langle\mathcal L_Kv,\phi_v\rangle\ge0.
\end{equation}
\end{lemma}
\begin{proof}
For fixed $\varepsilon>0$, the two-variable functions defining $\phi_u,\phi_v$ are Lipschitz on $[0,L]^2$ and vanish at $(0,0)$, so Lemma \ref{composition} gives admissibility. Here $u_\varepsilon,v_\varepsilon$ themselves do not belong to $\W$, but their differences coincide with those of $u,v$ respectively.

The four-case comparison below is a normalized version of the pointwise mechanism used in \cite[Proposition 3.1]{IannizzottoMugnai2025}; in order to avoid invoking boundary estimates for general kernels, we prove it directly. Fix $x,y$ and write
$a=u_\varepsilon(x)$, $b=u_\varepsilon(y)$, $c=v_\varepsilon(x)$, $d=v_\varepsilon(y)$, all positive, and set
$z_x=(a^p-c^p)^+$, $z_y=(b^p-d^p)^+$. We must show
\[
 D(x,y)=J_p(a-b)\left(\frac{z_x}{a^{p-1}}-\frac{z_y}{b^{p-1}}\right)
 -J_p(c-d)\left(\frac{z_x}{c^{p-1}}-\frac{z_y}{d^{p-1}}\right)\ge0
\]
where $J_p(t)=|t|^{p-2}t$.
If $a>c$ and $b>d$, then expanding gives
\begin{align*}
 D(x,y)={}&|a-b|^p+|c-d|^p
 -J_p(a-b)\left(\frac{c^p}{a^{p-1}}-\frac{d^p}{b^{p-1}}\right)\\
 &-J_p(c-d)\left(\frac{a^p}{c^{p-1}}-\frac{b^p}{d^{p-1}}\right)\ge0,
\end{align*}
the last step using \eqref{piconeeq} on the two cross terms separately. If $a>c$ and $b\le d$, then $z_y=0$ and
\[
 D(x,y)=z_x\left[J_p(1-b/a)-J_p(1-d/c)\right]\ge0,
\]
because $b/a\le d/c$ and $J_p$ is nondecreasing. If $a\le c$ and $b>d$, interchanging $x,y$ reduces to the previous case; if $a\le c$ and $b\le d$, then $D(x,y)=0$. Multiplying $D(x,y)\ge0$ by $K(x,y)$ and integrating gives \eqref{dsineq}. Finiteness of the integrals is ensured by the admissibility of the test functions together with H\"older's inequality.
\end{proof}

For each fixed  $f\in L^\infty(\Omega)$,  recall that the energy functional corresponding to~\eqref{eq:main} is
\begin{equation}\label{eq:energy}
E_f(u)= E(f,u)=\frac{1}{p}\Big([u]_K^p- \lambda\int_\Omega V(x)|u|^p\,dx\Big)+ \int_\Omega H(x,u)\,dx - \int_\Omega f u\,dx.
\end{equation}
Now, we give the proof of  Theorem~ \textup{\ref{thm:exist}}.
\begin{proof}[{\bf{Proof of Theorem~\textup{\ref{thm:exist}}}}]

By Lemma \ref{lem:ge0},  $H(x,t)= \int_0^t h(x,\tau)\,d\tau\ge 0$. It follows from (\ref{eq:energy}) and $0\le f\le M$ a.e.\ in $\Omega$ that
\begin{align*}
E_f(u) &\ge \frac{1}{p}\left( \iint_{\mathbb{R}^N\times\mathbb{R}^N} |u(x)-u(y)|^p K(x,y)\,dx\,dy - \lambda\int_\Omega V(x)|u|^p\,dx \right) - \int_\Omega f u\,dx\\
&\ge \frac{\delta^p}{p}[u]_K^p - M\|u\|_{L^1(\Omega)}\\
&\ge \frac{\delta^p}{p}[u]_K^p - C[u]_K.
\end{align*}
Since $p>1$, $E_f(u)\to\infty$ as $[u]_K\to\infty$.

The term $u\mapsto\frac{1}{p}[u]_K^p$ is convex and strongly continuous, hence weakly lower semicontinuous.
The term $u\mapsto- \lambda\int_\Omega V(x)|u|^p\,dx$ is weakly continuous by (V1)--(V3) and the compact embedding $W_0^{s,p}(\Omega)\hookrightarrow L^p(\Omega)$.
The term $u\mapsto\int_\Omega H(x,u)\,dx$ is weakly continuous by (H1) and (H3) and the compact embedding $W_0^{s,p}(\Omega)\hookrightarrow L^q(\Omega)$.
The term $u\mapsto\int_\Omega f u\,dx$ is linear and continuous.
Thus $E_f(u)$ is weakly lower semicontinuous.

Since $W_0^{s,p}(\Omega)$ is reflexive, $E_f(u)$ admits a global minimizer $u_f\in W_0^{s,p}(\Omega)$, which satisfies $ E^{'}_{f}(u_f)=0$ in $W^{-s,p'}(\Omega)$ and is therefore a weak solution of~\eqref{eq:main} in the sense of Definition~\ref{def:weak} (see Remark~\ref{rem:weak-well-posed}). Lemma \ref{bounded}, whose hypotheses are implied by those of the present theorem together with $0<f\le M$ a.e.\ in $\Omega$, then gives nonnegativity and boundedness.

We now prove uniqueness. Let $u_f,v_f$ be any two weak solutions of \eqref{eq:main}.  For simplicity, denote $u=u_f, v=v_f$    and set
$L=\max\{\|u\|_\infty,\|v\|_\infty\}$ and $\Omega_+=\{x\in\Omega:u(x)>v(x)\}$ in the following.
Use $u_\varepsilon,v_\varepsilon,z,\phi_u,\phi_v$ from Lemma \ref{ds}. Note that $\{z>0\}=\Omega_+$ is independent of $\varepsilon$. Testing the two equations with $\phi_u$ and $\phi_v$ respectively and subtracting, we obtain
\begin{equation}\label{ABC}
 0\le\langle\mathcal L_Ku,\phi_u\rangle-
 \langle\mathcal L_Kv,\phi_v\rangle=A_\varepsilon-B_\varepsilon-C_\varepsilon,
\end{equation}
where
\begin{align*}
 A_\varepsilon&=\lambda\int_{\Omega_+} V
 \left[\left(\frac{u}{u+\varepsilon}\right)^{p-1}
 -\left(\frac{v}{v+\varepsilon}\right)^{p-1}\right]z\,dx,\\
 B_\varepsilon&=\int_{\Omega_+}
 \left[\frac{h(x,u)}{(u+\varepsilon)^{p-1}}
 -\frac{h(x,v)}{(v+\varepsilon)^{p-1}}\right]z\,dx,\\
 C_\varepsilon&=\int_{\Omega_+} f
 \left[(v+\varepsilon)^{1-p}-(u+\varepsilon)^{1-p}\right]z\,dx.
\end{align*}
These three integrals are finite for each fixed $\varepsilon>0$. Clearly $A_\varepsilon,C_\varepsilon\ge0$. For $t>0$,
\[
 \frac{h(x,t)}{(t+\varepsilon)^{p-1}}
 =\frac{h(x,t)}{t^{p-1}}\left(\frac{t}{t+\varepsilon}\right)^{p-1}.
\]
Both factors on the right-hand side are nonnegative and nondecreasing, hence so is their product; at $t=0$ this quotient vanishes. Therefore (H2) and (H4) give $B_\varepsilon\ge0$. Consequently \eqref{ABC} implies
\begin{equation}\label{control}
 0\le C_\varepsilon\le A_\varepsilon.
\end{equation}

For $0<\varepsilon\le1$ the mean value theorem gives
\[
 0\le z\le p(L+1)^{p-1}(u-v)\quad\text{on }\Omega_+.
\]
The difference in brackets in $A_\varepsilon$ lies in $[0,1]$, so its integrand is dominated by the fixed $L^1$ function
\begin{equation}\label{domination}
 \lambda\|V\|_\infty p(L+1)^{p-1}(u-v)\mathbf1_{\Omega_+}.
\end{equation}
In particular, $A_\varepsilon$ and $C_\varepsilon$ are uniformly bounded.

We first rule out the possibility that $Z_{u,v}=\Omega_+\cap\{v=0\}$ has positive measure. On $Z_{u,v}$ we have $u>0$ and $f>0$, so the integrand of $C_\varepsilon$ is
\[
 f\left[\varepsilon^{1-p}-(u+\varepsilon)^{1-p}\right]
 \left[(u+\varepsilon)^p-\varepsilon^p\right]\longrightarrow+\infty.
\]
If $|Z_{u,v}|>0$, then applying Fatou's lemma to the nonnegative functions along any sequence $\varepsilon_j\downarrow0$ gives
$\liminf_j C_{\varepsilon_j}=+\infty$, contradicting \eqref{control} and the uniform boundedness. Hence $|Z_{u,v}|=0$.

Thus $u,v>0$ almost everywhere on $\Omega_+$, and the difference in brackets in $A_\varepsilon$ tends to zero. By \eqref{domination} and the dominated convergence theorem,
$A_\varepsilon\to0$; then \eqref{control} gives $C_\varepsilon\to0$. Applying Fatou's lemma once more to the nonnegative integrand of $C_\varepsilon$, we obtain
\[
 0\le\int_{\Omega_+} f\bigl(v^{1-p}-u^{1-p}\bigr)(u^p-v^p)\,dx
 \le\liminf_{\varepsilon\downarrow0}C_\varepsilon=0.
\]
On $\Omega_+$ this integrand is strictly positive almost everywhere, since $f>0$, $u>v>0$ and $p>1$. Hence $|\Omega_+|=0$, that is, $u\le v$ almost everywhere. Interchanging $u$ and $v$ yields $v\le u$ in the same way, so $u=v$.

Therefore, problem~\eqref{eq:main} admits a unique nonnegative weak solution $u_f\in W_0^{s,p}(\Omega)$  and $E(f,u_f)=\inf_{u\in W_0^{s,p}(\Omega)} E(f,u)$.
\end{proof}

\begin{remark}
The above proof is inspired by the $\varepsilon$-regularized quotient
method used in
\cite{MugnaiPinamontiVecchi,Kim}, but it employs different
comparison functions. The previous works used symmetric regularized
quotients together with upper truncations to control solutions that
might be unbounded. In the present paper, we first prove that the
solutions belong to $L^\infty(\Omega)$ and then set
\[
U=u+\varepsilon,
\qquad
W=v+\varepsilon,
\qquad
z_\varepsilon=(U^p-W^p)^+.
\]
We use
\[
\frac{z_\varepsilon}{U^{p-1}}
\qquad\text{and}\qquad
\frac{z_\varepsilon}{W^{p-1}}
\]
as test functions to obtain a one-sided comparison. The argument relies
only on a pointwise Picone-type inequality and the nonnegativity of
$K(x,y)$. It therefore applies to general measurable kernels that are
not necessarily translation invariant or radially symmetric. Moreover,
we only assume that
\[
t\longmapsto\frac{h(x,t)}{t^{p-1}}
\]
is nondecreasing. The strictness required for uniqueness is provided by
$f>0$, while the $p$-homogeneous potential term cancels in the comparison
argument. Therefore, the result in Theorem \ref{thm:exist}  is not a direct repetition of
the previous regularization arguments.
\end{remark}

%

\section{Solvability of the problem \eqref{eq:opt} and structure of the minimizer}\label{sec:prelim-lemmas}
{Suppose that $0<f_0\le M$ a.e.\ in $\Omega$. Then, for any $f\in\ov{\cR}\setminus\cR$, Lemma~\ref{lem:no-degenerate} gives $f>0$ a.e.\ in $\Omega$. By Theorem~\ref{thm:exist}, problem~\eqref{eq:main} admits a unique nonnegative weak solution $u_f\in W_0^{s,p}(\Omega)$ and $E(f,u_f)=\inf_{u\in W_0^{s,p}(\Omega)} E(f,u)$.} Therefore  
\begin{equation}\label{eq:psi}
\Psi(f)=\inf_{u\in W_0^{s,p}(\Omega)} E(f,u)
\end{equation}
is well defined.
Our main results in this section  are the following.

\begin{theorem}\label{thm:main}
Assume \textup{(K1)--(K2)}, \textup{(H1)--(H4)}, \textup{(V1)--(V3)} and $0<f_0\le M$ a.e.\ in $\Omega$.
Let $\cR$ be the rearrangement class of $f_0$ and $\Psi$ be defined by \eqref{eq:psi}. Then:
\begin{itemize}
\item[\rm{(i)}] the problem \eqref{eq:opt} is solvable, i.e., there exists $\hat f\in\cR$ such that $\Psi(\hat f)\le \Psi(f)$ for all $f\in\cR$;
\item[\rm{(ii)}] for any $\hat f\in\cR$ as in ${\rm(i)}$ and any $f\in\ov{\cR}\setminus\cR$, one has $\Psi(\hat f)<\Psi(f)$;
\item[\rm{(iii)}] any minimizer  $\hat f$ of $\Psi$ on  $\cR$ can be written as $\hat f = \eta\circ u_{\hat f}$ for some nondecreasing function $\eta:\R\to\R$.
\end{itemize}
\end{theorem}
In order to prove Theorem \ref{thm:main}, we prove some lemmas first.
\begin{lemma}[Absence of degenerate elements]\label{lem:no-degenerate}
Let $0<f_0\le M, a.e. ~x\in \Omega$  and  $\cR$ be the rearrangement class of $f_0$. Then for each $f\in\ov{\cR}$,   $0<f\le M, a.e. ~x\in \Omega$, especially,  $|\{x\in\Omega:\ f(x)=0\}|=0$.
\end{lemma}
\begin{proof}
By the convexity and the weak$^*$
lower semicontinuity of $f\mapsto\int_\Omega\zeta(f)\,dx$ for every convex function 
$\zeta:\R\to\R$, together with the fact that $\ov{\cR}$ is the closed convex hull of
$\cR$ (\cite[Theorem 6]{B87}, see also \cite[Lemma 2.14]{B89}), for each  $f\in\ov{\cR}$, 
we have
\[
\int_\Omega \zeta(f)\,dx\le\int_\Omega\zeta(f_0)\,dx.
\]
Taking $\zeta(t)=(-t)_+$ and using $f_0>0$ a.e.\ in $\Omega$ we get $f\ge 0$ a.e.\ in
$\Omega$; taking $\zeta(t)=t$ and $\zeta(t)=-t$ we get
$\int_\Omega f\,dx=\int_\Omega f_0\,dx$; and taking $\zeta(t)=(t-M)_+$ we get
\[
\int_\Omega (f-M)_+\,dx\le\int_\Omega (f_0-M)_+\,dx=0 ,
\]
that is, $f\le M$ a.e.\ in $\Omega$.

It remains to prove that $f>0$ a.e.\ in $\Omega$. Put $Z_f=\{x\in\Omega:\ f(x)=0\}$ and
$\vartheta=|Z_f|$. For every $g\ge0$ and every $\theta>0$ one has the pointwise identity
\[
\min(g,\theta)=g-(g-\theta)_+ .
\]
Applying it to $g=f_0$ and to $g=f$, and using $\int_\Omega f\,dx=\int_\Omega f_0\,dx$
together with the case $\zeta(t)=(t-\theta)_+$ of the display above, we obtain
\[
\int_\Omega \min(f_0,\theta)\,dx
=\int_\Omega f_0\,dx-\int_\Omega (f_0-\theta)_+\,dx
\le\int_\Omega f\,dx-\int_\Omega (f-\theta)_+\,dx
=\int_\Omega \min(f,\theta)\,dx .
\]
On the other hand $\min(f,\theta)=0$ on $Z_f$, while $\min(f,\theta)\le\theta$ on
$\Omega\setminus Z_f$; hence
\[
\int_\Omega \min(f,\theta)\,dx\le\theta\,(|\Omega|-\vartheta).
\]
Combining the last two estimates and dividing by $\theta>0$,
\[
\int_\Omega \min(f_0/\theta,\,1)\,dx
=\frac{1}{\theta}\int_\Omega \min(f_0,\theta)\,dx
\le |\Omega|-\vartheta .
\]
Since $f_0>0$ a.e.\ in $\Omega$, the integrands $\min(f_0/\theta,1)$ are dominated by $1$
and converge a.e.\ to $1$ as $\theta\to0^+$; dominated convergence therefore gives
$|\Omega|\le|\Omega|-\vartheta$, i.e.\ $\vartheta\le0$. As $\vartheta=|Z_f|\ge0$, we conclude
$\vartheta=0$, that is, $f>0$ a.e.\ in $\Omega$.
\end{proof}
The uniqueness part of Theorem~\ref{thm:exist} will be used in the following sharp form.
\begin{lemma}[Strong supporting hyperplane inequality]\label{lem:hyperplane}
{Let  $0<f_0\le M$ a.e.\ in $\Omega$  and  $\ov f\in\ov{\cR}$ be a minimizer of $\Psi$ over $\ov{\cR}$. Set $\ov u=u_{\ov f}$, the unique nonnegative weak solution of problem~\eqref{eq:main} with $f$ replaced by $\ov f$.}
Then, for every $f\in\ov{\cR}$,
\begin{equation}\label{eq:hyperplane}
\Psi(f)\ \le\ \Psi(\ov f)-\int_\Omega (f-\ov f)\,\ov u\,dx .
\end{equation}
Moreover, equality holds in \eqref{eq:hyperplane} if and only if $f=\ov f$ a.e.\ in $\Omega$.
\end{lemma}

\begin{proof}
 By Lemma~\ref{lem:no-degenerate}, {every $f\in\ov{\cR}$ satisfies} $0<f\le M$ a.e.\ in
$\Omega$, so $u_f$ is well defined and $\ov u$ is an admissible test function for
$E(f,\cdot)$. Hence, using $\Psi(f)=\inf_u E(f,u)$,
\[
\Psi(f)\ \le\ E(f,\ov u)
= E(\ov f,\ov u)-\int_\Omega(f-\ov f)\,\ov u\,dx
= \Psi(\ov f)-\int_\Omega(f-\ov f)\,\ov u\,dx ,
\]
where the last equality is the definition $\Psi(\ov f)=E(\ov f,u_{\ov f})$ of $\Psi$,
combined with the uniqueness of the minimizer of $E(\ov f,\cdot)$~ (i.e.,~$E_{\ov f}(\cdot) $) from
Theorem~\ref{thm:exist}.

Conversely, suppose that equality holds in \eqref{eq:hyperplane}. Then
$E(f,\ov u)=\Psi(f)$, so $\ov u$ minimizes $E(f,\cdot)$; by the uniqueness in
Theorem~\ref{thm:exist} we get $\ov u=u_f$. Testing the two equations $\eqref{eq:main}$
with sources $f$ and $\ov f$ at the common function $\ov u$ and subtracting, we obtain
\[
\int_\Omega (f-\ov f)\varphi\,dx=0\qquad\text{for every }\varphi\in W_0^{s,p}(\Omega),
\]
whence $f=\ov f$ a.e.\ in $\Omega$.
\end{proof}
\begin{lemma}[Weak$^*$ continuity]\label{lem:cont}
The functional $\Psi$ is sequentially weak$^*$ continuous on $\ov{\cR}$, where $\Psi$ is defined by \eqref{eq:psi}.
\end{lemma}

\begin{proof}
Let $(f_n)$ be a sequence in $\ov{\cR}$ such that $f_n \overset{*}{\rightharpoonup} f$ in $L^\infty(\Omega)$.
Denote $u_n = u_{f_n}$ and $u = u_f$.

\medskip\noindent\textit{Step 1: Upper bound.}
By definition of $\Psi$,
\[
\Psi(f_n) = \inf_{v\in W_0^{s,p}(\Omega)} E(f_n, v) \le E(f_n, u).
\]
Expanding the right-hand side,
\begin{align*}
E(f_n, u) &= \frac{1}{p}Q(u)+\int_\Omega \big[  H(x,u)-f_n u  \big]\,dx  \\
&= \frac{1}{p}Q(u)+\int_\Omega \big[ H(x,u)-f u   \big]\,dx   + \int_\Omega (f - f_n) u\,dx\\
&= \Psi(f) + \int_\Omega (f - f_n) u\,dx.
\end{align*}
Thus
\[
\Psi(f_n) \le \Psi(f) + \int_\Omega (f - f_n) u\,dx.
\]
Since $u\in L^1(\Omega)$ and $f_n \overset{*}{\rightharpoonup} f$, we have
\[
\lim_{n\to\infty} \int_\Omega (f_n - f) u\,dx = 0.
\]
Therefore
\begin{equation}\label{eq:limsup}
\limsup_{n\to\infty} \Psi(f_n) \le \Psi(f).
\end{equation}
This estimate is valid for the full sequence $(f_n)$.

\medskip\noindent\textit{Step 2: Boundedness of $\{u_n\}$.}
Each $u_n$ satisfies the weak formulation~\eqref{eq:weak}, namely
\[
\langle\cL_K u_n, \phi\rangle
- \lambda\int_\Omega V(x)|u_n|^{p-2}u_n \phi\,dx
+ \int_\Omega h(x,u_n)\phi\,dx
= \int_\Omega f_n \phi\,dx
\qquad \forall \phi\in W_0^{s,p}(\Omega).
\]
Taking $\phi = u_n$ and using the notation $Q(u_n)=[u_n]_K^p-\lambda\int_\Omega V|u_n|^p$,
\begin{equation}\label{eq:energy-norm}
Q(u_n) + \int_\Omega h(x,u_n)u_n\,dx = \int_\Omega f_n u_n\,dx.
\end{equation}

By (H2), $h(x,\cdot)$ is nondecreasing and $h(x,0)=0$, so $h(x,t)t \ge 0$ for all $t\in\R$.
Thus the second term on the left-hand side of \eqref{eq:energy-norm} is nonnegative.

For the right-hand side, $0 \le f_n \le M$ implies
\[
\int_\Omega f_n u_n\,dx \le M\|u_n\|_1.
\]
By the continuous embeddings $W_0^{s,p}(\Omega) \hookrightarrow L^p(\Omega) \hookrightarrow L^1(\Omega)$
(recall $|\Omega|<\infty$), there exists $C>0$ such that
\[
\|v\|_1 \le C[v]_K \qquad \forall v\in W_0^{s,p}(\Omega).
\]

From the lower bound $Q(u)^{1/p}\ge\delta[u]_K$ in Lemma~\ref{lem:phi-bound}, we have $[u_n]_K\le\delta^{-1}Q(u_n)^{1/p}$.
Dropping the nonnegative term in \eqref{eq:energy-norm} gives
\[
Q(u_n) \le \int_\Omega f_n u_n\,dx \le MC[u_n]_K \le MC\delta^{-1}Q(u_n)^{1/p}.
\]
Hence
\[
Q(u_n)^{(p-1)/p} \le MC\delta^{-1}.
\]
Since $p>1$, this implies $\{Q(u_n)\}$ is bounded, and by $Q(u_n)^{1/p}\ge\delta[u_n]_K$,
$\{[u_n]_K\}$ is bounded. Thus $\{u_n\}$ is bounded in $W_0^{s,p}(\Omega)$.

\medskip\noindent\textit{Step 3: Lower bound along a subsequence.}
Since $\{u_n\}$ is bounded in the reflexive space $W_0^{s,p}(\Omega)$ and $q<p_s^*$, after
passing to a further subsequence (still denoted by $\{u_n\}$) there exists $v\in W_0^{s,p}(\Omega)$ such that
\[
u_n \rightharpoonup v \ \text{ in } W_0^{s,p}(\Omega),\qquad
u_n \to v \ \text{ in } L^p(\Omega)\ \text{ and in } L^q(\Omega),\qquad
u_n(x) \to v(x) \ \text{a.e. in } \Omega.
\]

The functional $u \mapsto [u]_K^p$ is convex and continuous, hence weakly lower semicontinuous;
moreover $u_n \to v$ in $L^p(\Omega)$ and $V\in L^\infty(\Omega)$ give $\int_\Omega V|u_n|^p\,dx\to\int_\Omega V|v|^p\,dx$. Therefore
\[
\liminf_{n\to\infty} Q(u_n)
= \liminf_{n\to\infty}\Big( [u_n]_K^p - \lambda\int_\Omega V|u_n|^p \Big)
\ge [v]_K^p - \lambda\int_\Omega V|v|^p = Q(v).
\]

For the term involving $H$, using (H3) and the strong convergence $u_n \to v$ in $L^q(\Omega)$,
the Vitali convergence theorem yields
\[
\lim_{n\to\infty} \int_\Omega H(x,u_n)\,dx = \int_\Omega H(x,v)\,dx.
\]

For the linear term, since $f_n \overset{*}{\rightharpoonup} f$ and $u_n \to v$ strongly in $L^1(\Omega)$,
\begin{align*}
\lim_{n\to\infty} \int_\Omega f_n u_n\,dx
&= \lim_{n\to\infty} \left[ \int_\Omega f_n (u_n - v)\,dx + \int_\Omega f_n v\,dx \right]\\
&= 0 + \int_\Omega f v\,dx,
\end{align*}
because $|\int_\Omega f_n (u_n-v)\,dx| \le M\|u_n-v\|_1 \to 0$.

Now,
\begin{align*}
\liminf_{n\to\infty} \Psi(f_n)
&= \liminf_{n\to\infty} \left(\frac{1}{p}Q(u_n)+\int_\Omega H(x,u_n)\,dx- \int_\Omega f_n u_n\,dx  \right)\\
&\ge \liminf_{n\to\infty} \frac{1}{p}Q(u_n)+\lim_{n\to\infty} \int_\Omega H(x,u_n)\,dx-\lim_{n\to\infty} \int_\Omega f_n u_n\,dx \\
&\ge \frac{1}{p}Q(v)+\int_\Omega H(x,v)\,dx-\int_\Omega f v\,dx  \\
&= E(f,v).
\end{align*}

Since $\Psi(f) = \inf_u E(f,u)$, we have $E(f,v) \ge \Psi(f)$. Therefore, along the
subsequence selected in Step 3, $\liminf_{n\to\infty} \Psi(f_n) \ge \Psi(f)$.
The argument of Steps 2 and 3 uses only the uniform bound $0\le f_n\le M$ and the
convergence $f_n\overset{*}{\rightharpoonup} f$, both of which are inherited by every
subsequence of $(f_n)$; hence it can be repeated for any subsequence of $(f_n)$. Let us
show that this upgrades the last estimate to the whole sequence. Let
$\ell_*=\liminf_{n\to\infty}\Psi(f_n)$ and choose a subsequence with $\Psi(f_{n_k})\to\ell_*$; by
the preceding remark, after passing to a further subsequence, we may also assume that the
estimate holds along it. Then $\ell_*=\liminf_{k\to\infty}\Psi(f_{n_k})\ge\Psi(f)$, that is,
\begin{equation}\label{eq:liminf}
\liminf_{n\to\infty} \Psi(f_n) \ge \Psi(f).
\end{equation}

\medskip\noindent\textit{Step 4: Conclusion.}
From \eqref{eq:liminf} and \eqref{eq:limsup},
\[
\Psi(f) \le \liminf_{n\to\infty} \Psi(f_n) \le \limsup_{n\to\infty} \Psi(f_n) \le \Psi(f).
\]
Hence $\lim_{n\to\infty} \Psi(f_n) = \Psi(f)$, proving the sequential weak$^*$ continuity of $\Psi$.
\end{proof}

Now, we give the proof of  Theorem \ref{thm:main}.
\begin{proof}[{\bf{Proof of Theorem~\textup{\ref{thm:main}}}}]
$\ov{\cR}$ is weak$^*$ compact and $\Psi$ is weak$^*$ continuous (Lemma~\ref{lem:cont}), so
$\Psi$ attains its minimum on $\ov{\cR}$. Let $\ov f$ be a minimizer  of $\Psi$ on  $\ov{\cR}$  and set
$\ov u=u_{\ov f}$.

\medskip\noindent
{{(i).}} By Lemma~\ref{lem:linear-opt}(i) there exists $\hat f\in\cR$
maximizing the linear functional $f\mapsto\int_\Omega f\ov u\,dx$ over $\ov{\cR}$; in
particular
\[
\int_\Omega\hat f\ov u\,dx\ \ge\ \int_\Omega\ov f\ov u\,dx .\]
Consequently,
\[-\int_\Omega(\hat f-\ov f)\ov u\,dx\ \le\ 0 .
\]
{Applying Lemma~\ref{lem:hyperplane},   we obtain 
$$
\Psi(\hat f)\ \le\ \Psi(\ov f)-\int_\Omega (\hat f-\ov f)\,\ov u\,dx .
$$ }
Thus $\Psi(\hat f)\le\Psi(\ov f)$, and since $\ov f$ minimizes $\Psi$ over $\ov{\cR}$,  we conclude that
\begin{equation}\label{eq:min-equal}
\Psi(\hat f)=\Psi(\ov f)=\min_{\ov{\cR}}\Psi .
\end{equation}
As $\hat f\in\cR$ while $\cR\subseteq\ov{\cR}$, the equality of the minima
$\min_{\cR}\Psi=\min_{\ov{\cR}}\Psi$ follows, together with the fact that $\hat f$ minimizes
$\Psi$ over $\cR$. This proves (i).

\medskip\noindent
{ {(ii).}} Assume that  $\hat f\in\cR$ is a minimizer of $\Psi$ on  $\cR$ as in ${\rm(i)}$, so that
$\Psi(\hat f)=\min_{\cR}\Psi=\min_{\ov{\cR}}\Psi$ by \eqref{eq:min-equal}. Let
$f\in\ov{\cR}\setminus\cR$ and suppose, by contradiction, that $\Psi(f)=\Psi(\hat f)$;
then $f$ is a minimizer of $\Psi$ over $\ov{\cR}$ as well. Set $u=u_f$. By
Lemma~\ref{lem:hyperplane}, applied with $\ov f$ replaced by $f$, for every
$g\in\ov{\cR}$
\[
\Psi(g)\ \le\ \Psi(f)-\int_\Omega(g-f)u\,dx .
\]
By Lemma~\ref{lem:linear-opt}(i), choose $g\in\cR$ maximizing
$h\mapsto\int_\Omega hu\,dx$ over $\ov{\cR}$; then $\int_\Omega(g-f)u\,dx\ge0$, and hence
\[
\min_{\ov{\cR}}\Psi\ \le\ \Psi(g)\ \le\ \Psi(f)-\int_\Omega(g-f)u\,dx\ \le\ \Psi(f)
=\min_{\ov{\cR}}\Psi .
\]
All inequalities are therefore equalities. In particular equality holds in
\eqref{eq:hyperplane} with $\ov f$ replaced by $f$, so the equality clause of
Lemma~\ref{lem:hyperplane} gives $g=f$. Hence $f\in\cR$, which contradicts
$f\in\ov{\cR}\setminus\cR$. Therefore $\Psi(f)>\Psi(\hat f)$, which proves (ii).

\medskip\noindent
{ {(iii).}} Let $\hat f$ be a minimizer of $\Psi$ on  $\cR$ and $\hat u=u_{\hat f}$. By
Lemma~\ref{lem:hyperplane} (applied with $\ov f=\hat f$), for every $f\in\ov{\cR}$ with
$f\neq\hat f$,
\[
\Psi(f)\ \le\ \Psi(\hat f)-\int_\Omega(f-\hat f)\hat u\,dx ,
\qquad\text{so that}\qquad
\int_\Omega (f-\hat f)\hat u\,dx\ <\ 0 ,
\]
the strict inequality coming from (ii) when $f\notin\cR$ and from the equality clause of
Lemma~\ref{lem:hyperplane} when $f\in\cR$, $f\neq\hat f$. Equivalently,
$\int_\Omega f\hat u\,dx<\int_\Omega\hat f\hat u\,dx$: thus $\hat f$ is the unique
maximizer of $f\mapsto\int_\Omega f\hat u\,dx$ over $\ov{\cR}$. By
Lemma~\ref{lem:linear-opt}(ii), there exists a nondecreasing function $\eta:\R\to\R$ such
that $\hat f=\eta\circ\hat u$ a.e.\ in $\Omega$.
\end{proof}

\section*{Acknowledgments}
This work was supported by the National Natural Science Foundation of China (NSFC) under Grant No.~12261075 awarded to Chong Qiu and Grant No.~12601604 awarded to Pengcheng Wu.

\section*{Conflict of interest}
The authors declare no conflicts of interest.

\end{document}